\documentclass[sn-mathphys]{sn-jnl}

\jyear{2021}%

\theoremstyle{thmstyleone}%
\theoremstyle{thmstyletwo}%

\theoremstyle{thmstylethree}%

\newcommand{\ds}{\displaystyle}
\usepackage{subfigure}

\begin{document}

\title[Model-free saddle dynamics]{A model-free shrinking-dimer saddle dynamics for finding saddle point and solution landscape}

\author[1]{\fnm{Lei} \sur{Zhang}}\email{zhangl@math.pku.edu.cn}
\author[2]{\fnm{Pingwen} \sur{Zhang}}\email{pzhang@pku.edu.cn}
\author[3]{\fnm{Xiangcheng} \sur{Zheng}}\email{zhengxch@outlook.com}

\affil[1]{\orgdiv{Beijing International Center for Mathematical Research, Center for Machine Learning Research, Center for Quantitative Biology}, \orgname{Peking University}, \orgaddress{\city{Beijing}, \postcode{100871},  \country{China}}}

\affil[2]{\orgdiv{School of Mathematical Sciences, Laboratory of Mathematics and Applied Mathematics}, \orgname{Peking University}, \orgaddress{\city{Beijing}, \postcode{100871}, \country{China}}}

\affil[3]{\orgdiv{School of Mathematics}, \orgname{Shandong University}, \orgaddress{ \city{Jinan}, \postcode{250100}, 
\country{China}}}

\abstract{We propose a model-free shrinking-dimer saddle dynamics for finding any-index saddle points and constructing the solution landscapes, in which the force in the standard saddle dynamics is replaced by a surrogate model trained by the Gassian process learning. By this means, the exact form of the model is no longer necessary such that the saddle dynamics could be implemented based only on some observations of the force. This data-driven approach not only avoids the modeling procedure that could be difficult or inaccurate, but also significantly reduces the number of queries of the force that may be expensive or time-consuming. We accordingly develop a sequential learning saddle dynamics algorithm to perform a sequence of local saddle dynamics, in which the queries of the training samples and the update or retraining of the surrogate force are performed online and around the latent trajectory in order to improve the accuracy of the surrogate model and the value of each sampling. Numerical experiments are performed to demonstrate the effectiveness and efficiency of the proposed algorithm. }



\maketitle

\section{Introduction}\label{sec1}
Finding saddle points on complex energy landscapes or dynamical systems provides substantial physical and chemical properties and is thus an important problem in various fields such as nucleation and phase transformations in solid and soft matter \cite{ZhangChe,ZhangCheDu} and transition rates in chemical reactions and computational biology \cite{Nie2020,wang2021modeling,Yin2020nucleation}. The saddle points could be classified by the Morse index characterized by the maximal dimension of a subspace on which the Hessian $H(x)$ is negative definite \cite{Milnor}.  In practice, the index-1 saddle point represents the transition state connecting two local minima according to the transition state theory \cite{Mehta,ZhaRen}. The index-2 saddle points are particularly useful for providing trajectory information of chemical reactions in chemical systems \cite{Heidrich1986}. The excited states in quantum systems can also be characterized as saddle point configurations \cite{foresman1992toward}. 

There exist several algorithms of finding index-1 saddle points, e.g. \cite{EZho,Gou,ZhaDu}. However, the computation of high-index (index$>1$) saddle point is more difficult as it has multiple unstable eigen-directions and receives less attention, though the number of high-index saddle points are much larger than the number of local minima and index-1 saddle points on the energy landscape \cite{Chen2004,Li2001}. High-index shrinking-dimer saddle dynamics proposed in \cite{YinSISC}, which involves the force calculation (i.e. the negative gradient of the energy function), provides an effective means for finding any-index saddle points. This method is then combined with the downward search algorithm \cite{YinPRL,YinSCM} to construct the solution landscapes of both energy systems and dynamical (non-gradient) systems \cite{HanYin,han2021a,Luo2022, shi2022}. However, it is possible that the exact form of the force  is not given a priori in some cases such that we need to either investigate the modeling of the underlying processes or perform experiments in order to obtain the inquired values of the force in saddle dynamics. In practice, modeling complex problems may be difficult or inaccurate, which often restricts the computation of the saddle points and its applications. 

To resolve this issue, we employ a data-driven approach to propose a model-free saddle dynamics, in which the force in the original saddle dynamics is replaced by a surrogate model trained by the Gassian process learning. In the past few decades, the Gaussian process has been widely employed in extensive applications for constructing the surrogate models from the training data \cite{RaiPerKar,RasWil}. In particular, there exist some recent works on combining the Gaussian process with searching algorithms of index-1 saddle points \cite{Den,GuWan,KoiDag}. 
Here we adopt this idea in the computation of high-index saddle points. 
By this means, the exact form of the model is no longer necessary such that the saddle dynamics could be implemented based only on some observations of the force. Furthermore, the number of the force calculations during training is generally smaller than that in performing the saddle dynamics. 

Based on the proposed model-free saddle dynamics and its local nature, we adopt the sequential learning framework \cite{Pow} to develop a sequential learning saddle dynamics algorithm in which the queries of the training samples and the update or retraining of the surrogate force are performed online and around the latent trajectory. In this way, the accuracy of the surrogate model and the value of each sampling could be improved such that the queries of the force could be further reduced. The proposed method could be further combined with the downward search algorithm for construction of the model-free solution landscape \cite{YinPRL,YinSCM}, which provides a pathway map consisting of both saddle points and minima and avoids the sampling on energy landscape of the model system. Numerical experiments are performed to demonstrate the effectiveness and efficiency of the proposed algorithm in comparison with the standard saddle dynamics.

The rest of the paper is organized as follows: In Section 2 we propose a model-free shrinking-dimer saddle dynamics by incorporating the Gaussian process learning with the original saddle dynamics. In Section 3 we present a sequential learning saddle dynamics algorithm to compute high-index saddle points and construct the solution landscapes. Numerical experiments are performed in Section 4 and we draw concluding remarks in the last section.

\section{Model-free saddle dynamics}
We present a model-free shrinking-dimer saddle dynamics for finding high-index saddle points of complex energy functions or dynamical systems, in which the exact form of the force may not be given a priori  and is recovered from its observations at discrete locations. By this means, the saddle dynamics could be implemented based on these observations from, e.g., experiments or simulations, even without  knowing the formulations of the energy functions or dynamics. This data-driven approach not only avoids the modeling procedure that could be difficult or inaccurate, but could significantly reduce the number of queries of the true force that may be expensive or time-consuming in practical problems.  
\subsection{Saddle dynamics and its implementation}
We begin by introducing the high-index saddle dynamics proposed in \cite{YinSISC} to find an index-$k$ ($1\leq k\in\mathbb N$) saddle point of an energy function $E(x)$ 
\begin{equation}\label{Sadk}
\left\{
\begin{array}{l}
\ds \frac{dx}{dt} =\beta\bigg(I -2\sum_{j=1}^k v_jv_j^\top \bigg)F(x),\\[0.075in]
\ds \frac{dv_i}{dt}=\gamma \bigg( I-v_iv_i^\top-2\sum_{j=1}^{i-1}v_jv_j^\top\bigg)H(x)v_i,~~1\leq i\leq k.
\end{array}
\right.
\end{equation}
Here the force $F:\mathbb R^N\rightarrow \mathbb R^N$ is generated from the energy $E(x)$ by $F(x)=-\nabla E(x)$, $H(x):=-\nabla^2 E(x)$ corresponds to the Hessian of $E(x)$, $\beta$, $\gamma>0$ are relaxation parameters, $x$ represents the state variable and direction variables $\{v_i\}_{i=1}^k$ form a basis for the unstable subspace of the Hessian at $x$. 

More generally, (\ref{Sadk}) could be applied to find the saddle points of the dynamical systems $\dot{x}=F(x)$ which could be non-gradient, i.e., $F(x)$ is not a negative gradient of some energy function $E(x)$.  In this case, the Hessian in (\ref{Sadk}) should be replaced by the Jacobian of $F(x)$ or its symmetrization \cite{YinSCM,Z3}. Without loss of generality, we focus on the saddle dynamics of gradient systems in this paper.

As it is often expensive to calculate and store the Hessian, the dimer method \cite{HenJon,ZhaDu, ZhaSISC} is applied to approximate the multiplication of the Hessian and the vector $v$ as follows
\begin{equation}\label{Hzz}
H(x)v\approx \hat H(x,v,l):=\frac{F(x+lv)-F(x-lv)}{2l}
\end{equation}
where $2l$ refers to the dimer length for some $l>0$. 
Invoking this dimer approximation in (\ref{Sadk}) leads to the shrinking-dimer saddle dynamics \cite{YinSISC} 
\begin{equation}\label{sadk}
\left\{
\begin{array}{l}
\ds \frac{dx}{dt} =\beta\bigg(I -2\sum_{j=1}^k v_jv_j^\top \bigg)F(x),\\[0.075in]
\ds \frac{dv_i}{dt}=\gamma \bigg( I-v_iv_i^\top-2\sum_{j=1}^{i-1}v_jv_j^\top\bigg)\hat H(x,v_i,l),~~1\leq i\leq k,\\
\ds \frac{dl}{dt}=-\frac{dG}{dl}.
\end{array}
\right.
\end{equation}
Here an auxiliary function $G(l)$ defined on $[0,\infty)$ with $l = 0$ being its global minimum is introduced to control the dimer length. Following \cite{ZhaDu}, we take $G(l) = l^2/2$ to get an exponential decay of the
dimer length, while other choices such as $G(l) = l^3$ lead to a more gradual polynomial decay \cite{ZhaDuJCP}.

It is proved in \cite{YinSCM} that a linearly stable stationary point $x^*$ of the saddle dynamics (\ref{sadk}) is an index-$k$ saddle point satisfying $F(x^*)=0$. For practical implementation, we follow \cite{YinSISC,YinSCM} to get the scheme of the numerical solutions $\{x_n\}$, $\{v_{i,n}\}_{i=1}^k$ and $\{l_n\}$ for $n\geq 1$
 \begin{equation}\label{FDsadk}
\left\{
\begin{array}{l}
\ds x_{n} =x_{n-1}+\tau\beta\bigg(I -2\sum_{j=1}^k v_{j,n-1}v_{j,n-1}^\top \bigg)F(x_{n-1}),\\[0.1in]
\ds \tilde v_{i,n}=v_{i,n-1}+\tau\gamma\hat H(x_{n-1},v_{i,n-1},l_{n-1}),~~1\leq i\leq k,\\[0.2in]
\ds \{v_{i,n}\}_{i=1}^k=\text{GramSchmidt}(\{\tilde v_{i,n}\}_{i=1}^k),
\end{array}
\right.
\end{equation}
equipped with the prescribed initial values $x_0$, $l_0$ and (orthonormal) $\{v_{i,0}\}_{i=1}^k$.
 Here $l_{n-1}=l(\tau (n-1))$ for $n\geq 1$ are determined by solving the equation of $l$ in (\ref{sadk}) analytically under suitable choice of $G$, and the Gram-Schmidt orthonormalization procedure is applied in the third equation of (\ref{FDsadk}) to ensure the computational accuracy \cite{YinSISC}.  The scheme (\ref{FDsadk}) has been extensively applied to compute the saddle points and to construct the solution landscapes \cite{HanYin,han2021a,YinPRL}. However, implementing this scheme requires the values of $F$ at $x_{n-1}$ and $x_{n-1}\pm lv_{i,n-1}$ for $1\leq i\leq k$ at each time step $t_n$, which may be difficult to evaluate if the model (e.g. the form of $F$) is unknown or the query of $F$ is expensive or time-consuming. 

\subsection{Surrogate force via Gaussian process}
To accommodate the concerns mentioned above, we intend to replace the true force $F(x)=[F_1(x),\cdots,F_N(x)]^\top$ in (\ref{FDsadk}) by a surrogate force $\mathcal F(x)=[\mathcal F_1(x),\cdots,\mathcal F_N(x)]^\top$, which is learned from the observations of $F$ at training locations via the multi-output Gaussian process, see e.g. \cite{LiuCai}. Suppose  $\mathcal F(x)$ satisfies the following Gaussian process
\begin{equation}\label{calF}
\mathcal F(x)\sim \mathcal{GP}(0,\mathcal K(x,x')) 
\end{equation}
where the multi-output covariance $\mathcal K(x,x')\in \mathbb R^{N\times N}$ is defined as
$$\mathcal K(x,x')=
\begin{bmatrix}
		k_{1,1}(x,x')&\cdots & k_{1,N}(x,x')\\
		 \vdots &\ddots&\vdots\\
		 k_{N,1}(x,x')&\cdots &k_{N,N}(x,x')
	\end{bmatrix}
 $$
with the element $k_{p,q}$ ($1\leq p,q\leq N$) corresponding to the covariance between $\mathcal F_p(x)$ and $\mathcal F_q(x')$. Denote the hyper-parameters in $k_{p,q}$ ($1\leq p,q\leq N$) as $\theta$.

Given the $m$ observations $Y:=\{y^1,\cdots,y^m\}$ of $F$ at the training data set $X:=\{x^t_1,\cdots,x^t_m\}$ where $y^j=[y^j_1,\cdots,y^j_N]^\top$ for $1\leq j\leq m$ and the observation of the $i$th output $y_i$ ($1\leq i\leq N$) is assumed to satisfy $y_i(x)=\mathcal F_i(x)+\varepsilon_i$ with the independent and identically distributed Gaussian noise $\varepsilon_i\sim \mathcal{N}(0,\sigma_{s,i}^2)$. Then the hyper-parameters $\theta$ and the variances $\{\sigma_{s,i}\}_{i=1}^N$ could be inferred by the standard maximum likelihood estimation of maximizing the marginal likelihood  $p(Y\vert X;\theta,\{\sigma_{s,i}\}_{i=1}^N)$ \cite{LiuCai,LiuZhu,RasWil}.

After inferring these parameters, the posterior distribution of $\mathcal F$ at any test point $x_*$ could be analytically derived as \cite{LiuCai,RasWil}
$$\mathcal F(x_*)\vert X,Y,x_*\sim \mathcal{N}(\mu_*,\Sigma_*) $$
where the predicted mean and variance are presented as
$$ \mu_*=K_*^\top (K(X,X)+\Sigma)^{-1}\bar y$$
and 
$$\Sigma(x_*)=\mathcal K(x_*,x_*)-K_*^\top (K(X,X)+\Sigma)^{-1}K_* .$$
Here $\bar y\in\mathbb R^{mN\times 1}$ is a rearrangement of the observations in $Y$ defined by 
$$\bar y =[y^1_1,\cdots,y^m_1,y^1_2,\cdots,y^m_2,\cdots,y^1_N,\cdots,y^m_N]^\top,$$
 $K_*\in\mathbb R^{mN\times N}$ is defined by
 $$K_*=
\begin{bmatrix}
		K_{1,1}(X,x_*)&\cdots & K_{1,N}(X,x_*)\\
		 \vdots &\ddots&\vdots\\
		 K_{N,1}(X,x_*)&\cdots &K_{N,N}(X,x_*)
	\end{bmatrix}
	 $$
  where 
  $$K_{p,q}(X,x_*)=[k_{p,q}(x^t_1,x_*),\cdots,k_{p,q}(x^t_m,x_*)]^\top$$ for $1\leq p,q\leq N$, $\Sigma=\Sigma_s\bigotimes I_m\in\mathbb R^{mN\times mN}$ is the noise matrix with 
  $$\Sigma_s=
  \begin{bmatrix}
		\sigma_{s,1}^2& & \\
		  &\ddots&\\
		 & &\sigma_{s,N}^2
	\end{bmatrix},
	$$
   and the symmetric block matrix $K(X,X)\in\mathbb R^{mN\times mN}$ is given as
$$K(X,X)=
\begin{bmatrix}
		K_{1,1}(X,X)&\cdots & K_{1,N}(X,X)\\
		 \vdots &\ddots&\vdots\\
		 K_{N,1}(X,X)&\cdots &K_{N,N}(X,X)
	\end{bmatrix}
	 $$
with 
$$K_{p,q}(X,X)=
\begin{bmatrix}
		k_{p,q}(x^t_1,x^t_1)&\cdots & k_{p,q}(x^t_1,x^t_m)\\
		 \vdots &\ddots&\vdots\\
		 k_{p,q}(x^t_m,x^t_1)&\cdots &k_{p,q}(x^t_m,x^t_m)
	\end{bmatrix}
	$$
for $1\leq p,q\leq N$.

In practice, we choose $\mu_*$ as the value of $\mathcal F(x_*)$, i.e. we denote 
$\mathcal F(x_*)=\mu_*,$
 and thus the predicted value of $ F(x_*)$, and the diagonal entries of $\Sigma(x_*)$ represent its variances. Then we replace the true force $F(x)$ by its surrogate model $\mathcal F(x)$ in the saddle dynamics (\ref{sadk}) to get the following model-free saddle dynamics
\begin{equation}\label{modsadk}
\left\{
\begin{array}{l}
\ds \frac{dx}{dt} =\beta\bigg(I -2\sum_{j=1}^k v_jv_j^\top \bigg)\mathcal F(x),\\[0.075in]
\ds \frac{dv_i}{dt}=\gamma \bigg( I-v_iv_i^\top-2\sum_{j=1}^{i-1}v_jv_j^\top\bigg)\mathcal H(x,v_i,l),~~1\leq i\leq k,\\
\ds \frac{dl}{dt}=-\frac{dG}{dl}
\end{array}
\right.
\end{equation}
with 
$$
\mathcal H(x,v_i,l):=\frac{\mathcal F (x+lv_i)-\mathcal F(x-lv_i)}{2l}.
 $$
  The corresponding numerical scheme takes an analogous form as (\ref{FDsadk})
  \begin{equation}\label{modFDsadk}
\left\{
\begin{array}{l}
\ds x_{n} =x_{n-1}+\tau\beta\bigg(I -2\sum_{j=1}^k v_{j,n-1}v_{j,n-1}^\top \bigg)\mathcal F(x_{n-1}),\\[0.1in]
\ds \tilde v_{i,n}=v_{i,n-1}+\tau\gamma\mathcal H(x_{n-1},v_{i,n-1},l_{n-1}),~~1\leq i\leq k,\\[0.2in]
\ds \{v_{i,n}\}_{i=1}^k=\text{GramSchmidt}(\{\tilde v_{i,n}\}_{i=1}^k).
\end{array}
\right.
\end{equation}

\section{A sequential learning algorithm}
We propose a sequential leaning algorithm involving the model-free saddle dynamics scheme (\ref{modFDsadk}) in practical computations of saddle points and solution landscapes. The motivation is that due to the local nature of the saddle dynamics, the training data located far from the dynamical path of saddle dynamics may have less contributions or even introduce errors in predicting the force $F(x)$. Ideally, the training data should be sampled near the dynamical path, which, however, is not known a priori that makes the task of generating the training data for training the surrogate force nontrivial. To accommodate this issue, we adopt the sequential learning technique to perform the model-free saddle dynamics through an active learning framework, where the acquisition of training samples and the validation and update of the surrogate force are performed ``online'' (during optimization). By this means, we expect to reduce the number of training points while preserving the accuracy of the surrogate model and thus the computed saddle points.

The basic idea of the sequential learning is to divide the learning-based optimization into several trust region optimization subproblems. For the $j$th subproblem, we perform the Gaussian process learning to train the force $F(x)$ in a hypercube trust region 
$$Q(x^j_c,\Delta^j):=\{x\vert\|x-x^j_c\|_\infty\leq \Delta^j\}$$
 for the center $x^j_c$ and trust region length $\Delta^j$ based on a training data set with $N_{sam}$ samples $D:=\{(x^t_1,y_1),\cdots,(x^t_{N_{sam}},y_{N_{sam}})\}$  and then perform the surrogate-force-based saddle dynamics (\ref{modFDsadk}) inside this region for $n=1,2\cdots$ until the $x_{n}$ gets out of this region for some $n=n_*$. 
 
 Then we check the reliability of $\mathcal F$ at $x_{n_*}$ in order to determine the center $x^{j+1}_c$ and the trust region length $\Delta^{j+1}$ of the $(j+1)$th subproblem. We compute the maximum norm of the diagonal of the covariance matrix $\Sigma(x_{n_*})$ of $\mathcal F(x_{n_*})$, that is, 
 \begin{equation}\label{r}
 r:=\|\mbox{diag}(\Sigma(x_{n_*}))\|_\infty.
 \end{equation}
  Given the lower and upper tolerances $0<$tol$_l<$tol$_u$, we encounter three cases:
  \begin{itemize}
\item   If $r<\mbox{tol}_l$, which means that the surrogate model is fairly reliable, we could set $x^{j+1}_c$ as $x_{n_*}$ and enlarge the trust region length in the $(j+1)$th subproblem.  

\item For the case $r>\mbox{tol}_u$, then the surrogate model may not be reliable such that the $j$th subproblem should be solved again using the retrained surrogate force  under the updated training data set with $N_{new}$ newly added samples. In this case, we set $x^{j+1}_c$ back to $x^j_{c}$ and shrink the trust region length to improve the accuracy.  

\item For tol$_l\leq r\leq$tol$_u$, we set $x^{j+1}_c$ as $x_{n_*}$ and keep the  trust region length unchanged to solve the $(j+1)$th subproblem.
\end{itemize}

In all three cases, the training data obtained before will be inherited if they locate in the updated region, which helps to improve the accuracy of the surrogate model.  The algorithm will be terminated if the error of two adjacent numerical solutions of $x$ is smaller than the tolerance tol$_x$ or the total number of steps $N_{t}$ exceeds the prescribed value $N_m$.  We summarize this algorithm in \textbf{Algorithm 1}, which computes an index-$k$ saddle point and count the total number $N_f$ of queries of $F$.

\begin{algorithm}[H]
\caption{Sequential learning algorithm of finding an index-$k$ saddle point}\label{algo1}
{\small
\begin{algorithmic}[1]
\State \textbf{Step 1: Initialization}
\State Initial data: $x_0$, $\{v_{i,0}\}_{i=1}^k$, $l_0$, $\beta$, $\gamma$, $\tau$, tol$_l$, tol$_u$, $N_m$, $N_t=0$, $N_f=0$, tol$_x$, $N_{sam}$, $N_{new}$, $x_c$, $\Delta$.
\State Initial sampling: Sample $N_{sam}$ points $X=\{x^t_1,\cdots,x^t_{N_{sam}}\}$ in the trust region $Q(x_c,\Delta)$ and then evaluate the corresponding values of $F$ as $Y=\{F(x^t_1),\cdots,F(x^t_{N_{sam}}))\}$ to get the training date set $D=\{X,Y\}$. 
\State $N_f\Leftarrow N_f+N_{sam}$.

\State \textbf{Step 2: Implement model-free saddle dynamics within $Q(x_c,\Delta)$}
\State Infer the hyper-parameters $\theta$ and the variances $\{\sigma_{s,i}\}_{i=1}^N$ by maximum likelihood estimation with the training set $D$ to obtain $\mathcal F$.
\State Implement saddle dynamics (\ref{modFDsadk}) using $\mathcal F$ until one of the following cases occurs: \\
(\textbf{a}) $x_{n}\not\in Q(x_c,\Delta)$ for some $n=n_*$; (\textbf{b}) $N_t>N_m$; (\textbf{c}) $\|x_n-x_{n-1}\|_{\infty}\leq$tol$_x$ for some $n$. \\
For the cases (\textbf{b}) and (\textbf{c}), terminate the algorithm and output the latest numerical solution of $x$ and $N_f$.
\State \textbf{Step 3: Update of $x_c$, $\Delta$ and $D$}
\State Compute $r$ in (\ref{r}).
\If{$r<$tol$_l$}
       $x_c \Leftarrow  x_{n_*}$ and $\Delta \Leftarrow 2\Delta$;
\Else 
\If{$r>$tol$_u$}
       $\Delta \Leftarrow \Delta/2$;
  \Else\,
   $x_c \Leftarrow  x_{n_*}$.
  \EndIf
\EndIf
\State Resample $N_{new}$ training points in $Q(x_c,\Delta)$ as in \textbf{Step 1} to get a training data set $D_{new}$.
\State $N_f\Leftarrow N_f+N_{new}$,  $D \Leftarrow D\vert_{x\in Q(x_c,\Delta)}\cup D_{new}$ and jump to \textbf{Step 2}. 
\end{algorithmic}
}
\end{algorithm}

Based on \textbf{Algorithm 1}, we could apply the downward  or upward search algorithms proposed in \cite{YinPRL,YinSCM} for construction of solution landscapes. For instance, the downward search algorithm aims to search lower-index saddles and stable states from a high-index saddle point. Given an index-$k$ saddle point and its $k$ unstable directions denoted by $(x^*, v_1, \cdots, v_k)$ as a parent state, we could apply this algorithm to search for an index-$m$ saddle point with $0 \le m < k$. Typically, the initial value is chosen as a perturbation of the parent state $x^*$  along the direction $v_i $ for some $m < i \le k$ with the directions $\{v_j\}_{j=1}^m$ as the initial vectors to start the saddle dynamics. By repeating this procedure on the newly found saddle points, the landscape under the given index-$k$ saddle point could be constructed completely.


\section{Numerical experiments}
We present numerical examples to substantiate the effectiveness and the computational cost of the sequential-learning Gaussian process saddle dynamics (GPSD) proposed by \textbf{Algorithm 1} in finding saddle points, in comparison with the standard saddle dynamics (SD).  Here the computational cost is characterized by the number $N_f$ of the force evaluations in order to demonstrate that the GPSD is more efficient than the SD in practical problems when the query of the force is expensive or time-consuming. In numerical experiments, we follow the conventional treatment to take $k_{p,q}$ for $1\leq p,q\leq N$ as the squared exponential covariance function with the hyper-parameters $\theta:=(\sigma_f,\sigma_l)$
$$k(x,x'):=\sigma_f^2 \mbox{exp}\Big(-\frac{\|x-x'\|_2^2}{2\sigma_l^2}\Big) .$$
\subsection{A Rosenbrock type function}
We compute the saddle points of the following four-dimensional Rosenbrock type function
\begin{align*}
 &E(x_1,x_2,x_3,x_4)\\
 &\quad\quad\quad:=a(x_4-x_3^2)^2+b(x_3-x_2^2)^2+c(x_2-x_1^2)^2+d(1-x_1)^2.
 \end{align*}
We set $\tau=l_0=0.01$, tol$_l$=0.05, tol$_u$=0.15, $N_m=20000$, tol$_x$=$1\times 10^{-6}$, $N_{sam}=N_{new}=100$, $x_0=(0.7, 0.8, 1.2, 0.7)$, $\Delta=0.025$ and the dimer length with polynomial decay 
$$l_n=\frac{l_0}{1+(n\tau)^2}.$$ The initial values $\{v_{i,0}\}_{i=1}^k$ are chosen as the orthonormal eigenvectors of the Hessian of $E$ at $x_0$. The Latin hypercube sampling technique \cite{Tan} is applied to generate the training data locations. Under different coefficients $a$, $b$, $c$ and $d$, the index of the saddle point $(1,1,1,1)$ of $E$ is different, and we perform both methods for the four cases in Table \ref{table1}.

\begin{table}[h]
\setlength{\abovecaptionskip}{0pt}
\centering
\caption{Index of the saddle point $(1,1,1,1)$ for $E$ under different parameters.}
\begin{tabular}{cccccc} \cline{1-6}
&$a$& $b$ & $c$ & $d$ &Index\\ \cline{1-6}		
(i)&$-0.5$&	$0.5$&	$0.5$	&$2$	& 1\\
(ii)&$-0.5$&	$0.5$&	$-0.5$	&$2$	& 2\\
(iii)&$-0.5$&	$-0.5$&	$-0.5$	&$2$	& 3\\
(iv)&$-0.5$&	$-0.5$&	$-0.5$	&$-2$	& 4\\
				\hline
			\end{tabular}
			\label{table1}
		\end{table}
		
	Numerical results are presented in Tables \ref{table2}--\ref{table3}, which indicate that both GPSD and SD generate satisfactory approximations for the saddle point $(1,1,1,1)$ with different indexes, while the GPSD requires much less queries of the true force $F$, which demonstrate the potential of the GPSD 
in practical problems where the acquisition of the force is expensive or time-consuming.
	\begin{table}[h]
\setlength{\abovecaptionskip}{0pt}
\centering
\caption{Comparison of SD and GPSD for approximating the saddle point $(1,1,1,1)$ of $E$}
\begin{tabular}{cll} \cline{1-3}
&\hspace{0.56in} SD&\hspace{0.45in} GPSD  \\ \cline{1-3}		
(i)&(1.000,    0.999,    0.998,   0.995)&	$(1.022,1.016,1.018,1.032)$	\\
(ii)&(1.000,    1.0000,    1.000,    1.000) &(1.011, 0.997,    1.004,    1.006)	\\
(iii)&(1.000,   1.000,   1.000,    0.999) &(1.002,  0.991,   1.012,    0.999)   \\
(iv)&(1.000,   0.999,    0.998,    0.995)&(0.990,1.002,    0.993,    0.991)	\\
				\hline
			\end{tabular}
			\label{table2}
		\end{table}	
		
		\begin{table}[h]
\setlength{\abovecaptionskip}{0pt}
\centering
\caption{Comparison of the number of force evaluation $N_f$ in SD and GPSD for approximating the saddle point of $E$}
\begin{tabular}{ccc} \cline{1-3}
&$N_f$ in SD&$N_f$ in GPSD  \\ \cline{1-3}		
(i)&$50673$&5200	\\
(ii)&13995 &2700		\\
(iii)&23548 &3100	     \\
(iv)&155502&5200		\\
				\hline
			\end{tabular}
			\label{table3}
		\end{table}

\subsection{A modified codesign problem}		\label{sec42}
We compute the saddle point of the following modification of the codesign problem \cite{All,Des} that simultaneously optimizes the plant design variable $a$ and the control design variable $u(t)$ defined on $[0,0.1]$ in the plant design objective $J=J(a,u)$
\begin{equation}\label{J}
J(a,u):=-\frac{a^2}{2}+\frac{1}{2}\int_0^{0.1}\xi^2(t)u^2(t)dt. 
\end{equation}
Different from the most conventional problems, in which the energy or cost functionals are exactly given in priori, there is an unknown state variable $\xi(t)$ in $J$ that relates to both $a$ and $u(t)$ via the underlying mechanism. In other words, for each query of $J$ or the corresponding negative gradient $F=[a,-\xi(t)^2u(t)]^\top$ at some $(a,u(t))$ in each iteration of SD, we need to compute $\xi(t)$ by simulating or solving the governing process related to the current $(a,u(t))$, which is in general computationally expensive. Therefore, we expect to apply GPSD, which constructs a surrogate mapping from $(a,u(t))$ to $F$ to partly circumvent queries of $\xi(t)$ in evaluating $F$,  to reduce the number of simulating or computing the governing process during the SD iterations.

 From $F=[a,-\xi(t)^2u(t)]^\top$ we find that $x^*:=[0,0]^\top$ is an index-1 saddle point of $J$ and we will apply both SD and GPSD for its computation.  For illustration, we follow \cite[Example 1]{Des} to select the governing process as a system of nonlinear differential equations on $(0,0.1]$
\begin{equation}\label{sys}
\left\{
\begin{array}{l}
\ds \frac{d\eta}{dt} =-a\eta+\xi^2,\\[0.1in]
\ds \frac{d\xi}{dt}=\eta-2a^2\xi-\eta^2+u,
\end{array}
\right.
\end{equation}
equipped with the initial conditions $\eta(0)=\xi(0)=1$. In practice, much more complicated equations or processes than (\ref{sys}) could appear.  
We set $\tau=l_0=0.0025$, tol$_l$=0.05, tol$_u$=0.15, $N_m=20000$, tol$_x$=$1\times 10^{-6}$, $N_{sam}=N_{new}=300$, $\Delta=0.1$ and the dimer length with polynomial decay as before.
We discretize the domain $[0,0.1]$ of $u$, $\eta$ and $\xi$ by a uniform partition with mesh size $0.01$, and approximate their values on this partition where the explicit difference method is used to solve (\ref{sys}).
The initial value $v_{1,0}$ of saddle dynamics is chosen as the normalized vector of $[1,\cdots,1]^\top$, and different initial values are applied as follows
$$(\mbox{i}) \,x_0=[0.2,\cdots,0.2]^\top;~(\mbox{ii})\, x_0=[0.1,\cdots,0.1]^\top;~(\mbox{iii})\, x_0=[0.05,\cdots,0.05]^\top.$$  
The Latin hypercube sampling technique \cite{Tan} is applied to generate the training data locations. We measure the number $N_s$ of simulating (\ref{sys}) and the errors between $x^*$ and its output numerical approximation $x_F$ in both methods in Table \ref{table33}, which indicate that both methods generate satisfactory results while the GPSD performs much less simulations of the underlying process. Furthermore, as the initial value $x_0$ approaches $x^*$ from (i) to (iii), the $N_s$ of GPSD rapidly decreases while the $N_s$ of SD does not decrease evidently. All these observations demonstrate the advantages of the proposed GPSD method.

		\begin{table}[h]
\setlength{\abovecaptionskip}{0pt}
\centering
\caption{Comparison of the number $N_s$ of simulating the underlying process (\ref{sys}) and the errors $\|x_F-x^*\|_{\infty}$ in SD and GPSD for approximating the saddle point $x^*$ of $J$}
\begin{tabular}{ccccc} \cline{1-5}
& $N_s$ in SD&$\|x_F-x^*\|_{\infty}$ in SD& $N_s$ in GPSD&$\|x_F-x^*\|_{\infty}$ in GPSD \\ \cline{1-5}		
(i)&10494&3.98$\times 10^{-4}$&4800&8.22$\times 10^{-3}$\\
(ii)&9618 &3.98$\times 10^{-4}$&2700&9.47$\times 10^{-4}$\\
(iii)& 8775&3.99$\times 10^{-4}$&1200&6.84$\times 10^{-3}$	     \\

				\hline
			\end{tabular}
			\label{table33}
		\end{table}	
		
\subsection{A nonlocal phase-field model}
\begin{figure}[H]
	\setlength{\abovecaptionskip}{0pt}
	\centering	
	\subfigure[Initial value $u_0(x)$]{
	\includegraphics[width=1.8in,height=1.4in]{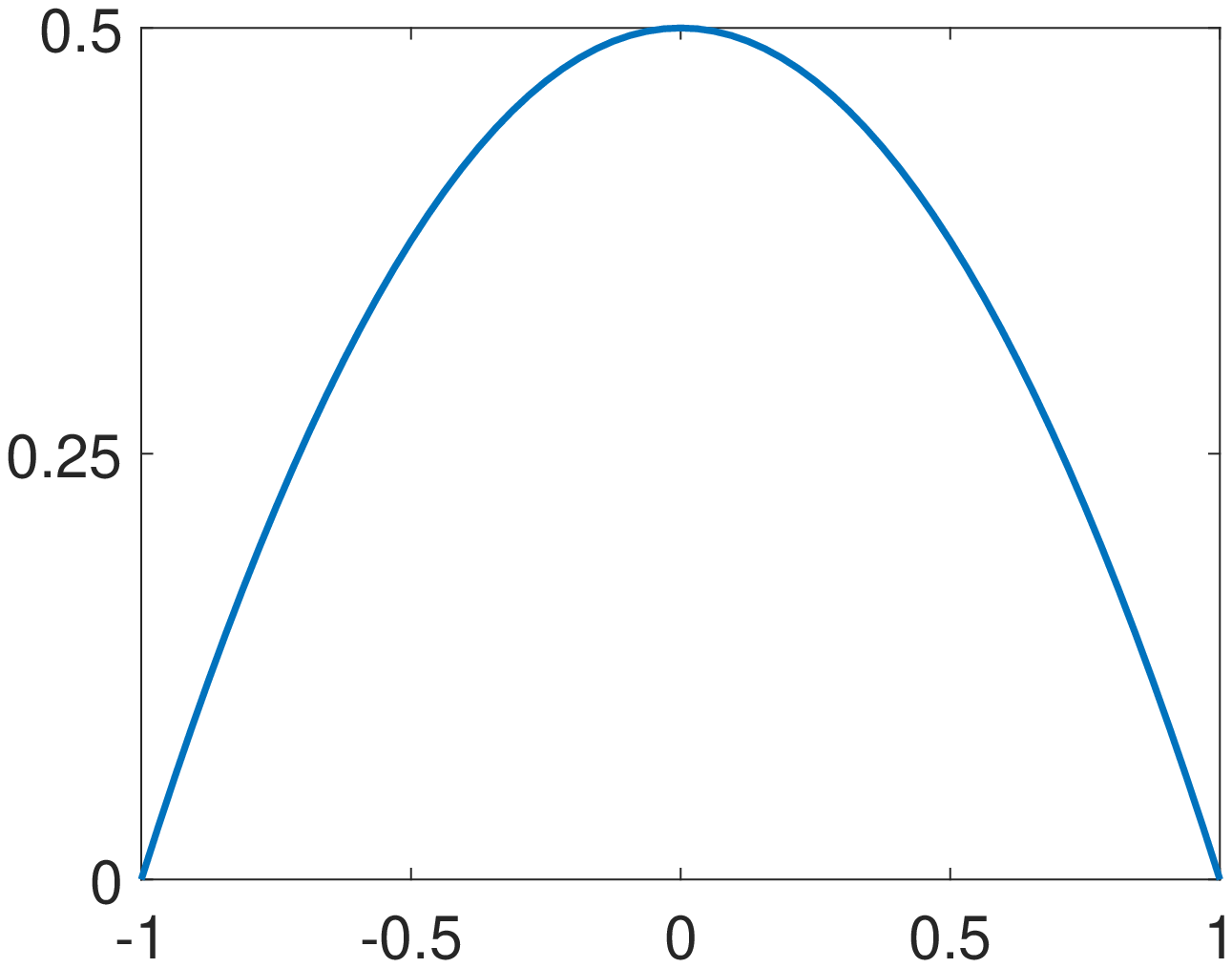}
	}
	\subfigure[An index-1 saddle point]{
	\includegraphics[width=1.8in,height=1.4in]{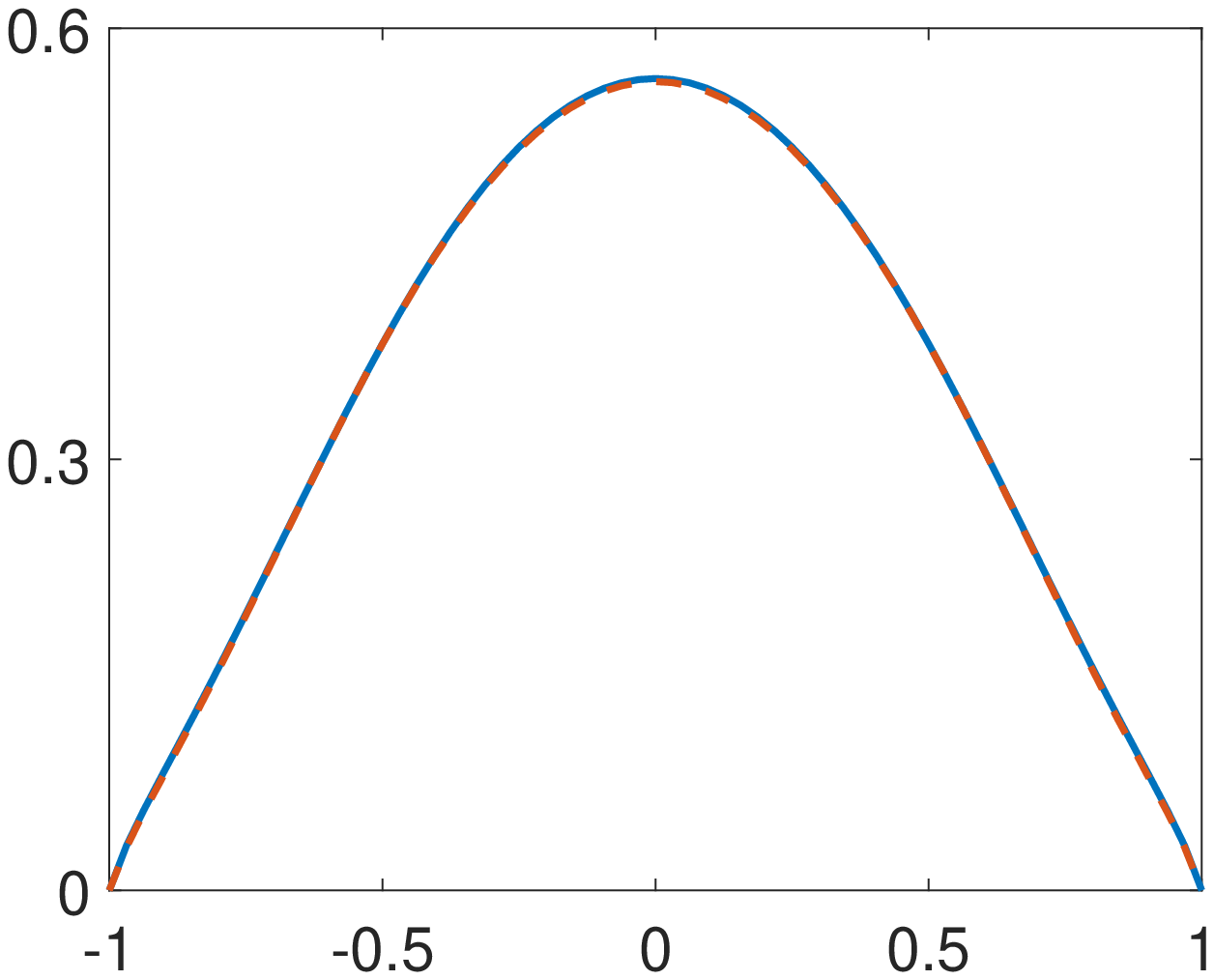}
	}\\
\subfigure[An index-2 saddle point]{\includegraphics[width=1.8in,height=1.4in]{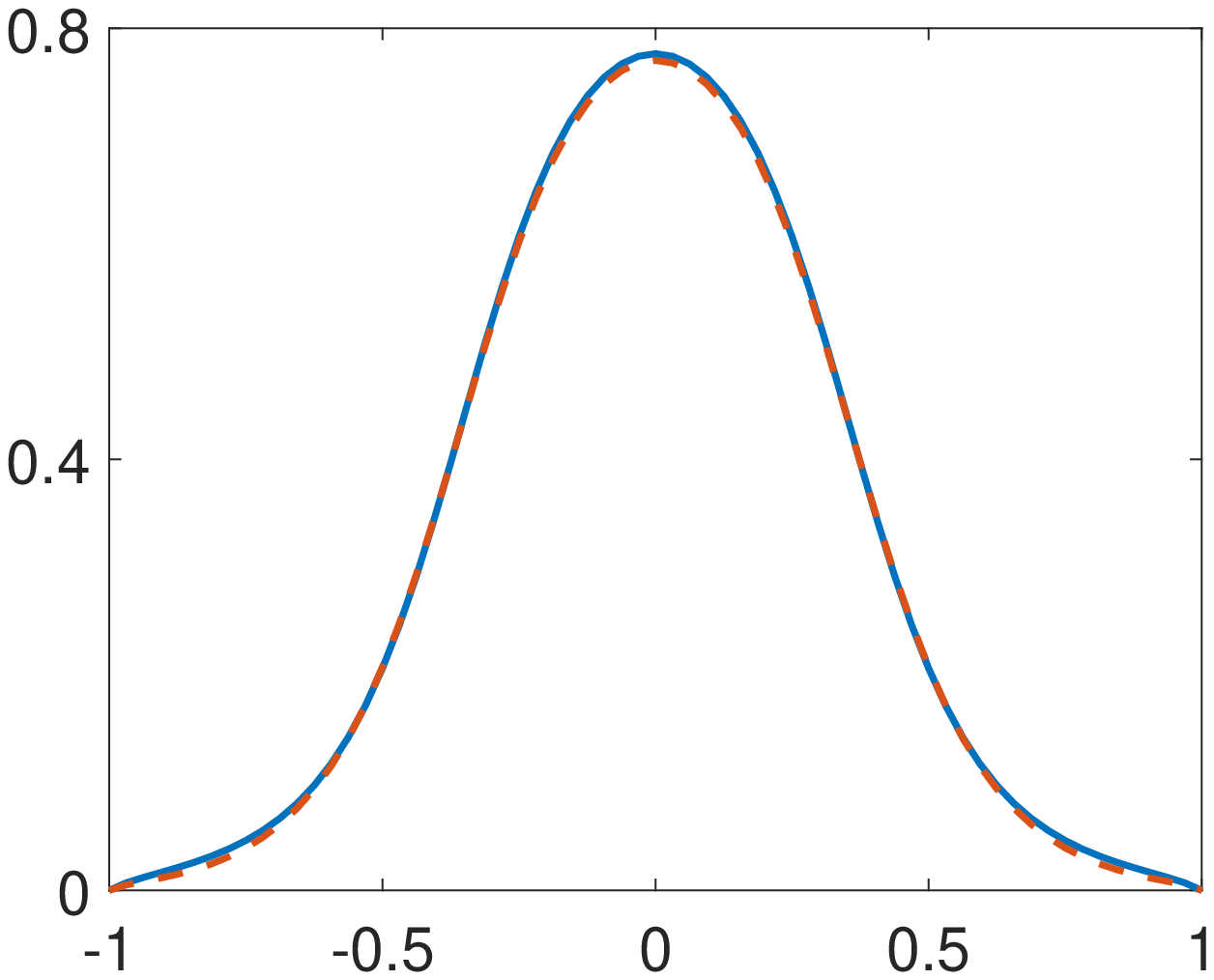}
}
\subfigure[An index-3 saddle point]{
\includegraphics[width=1.8in,height=1.4in]{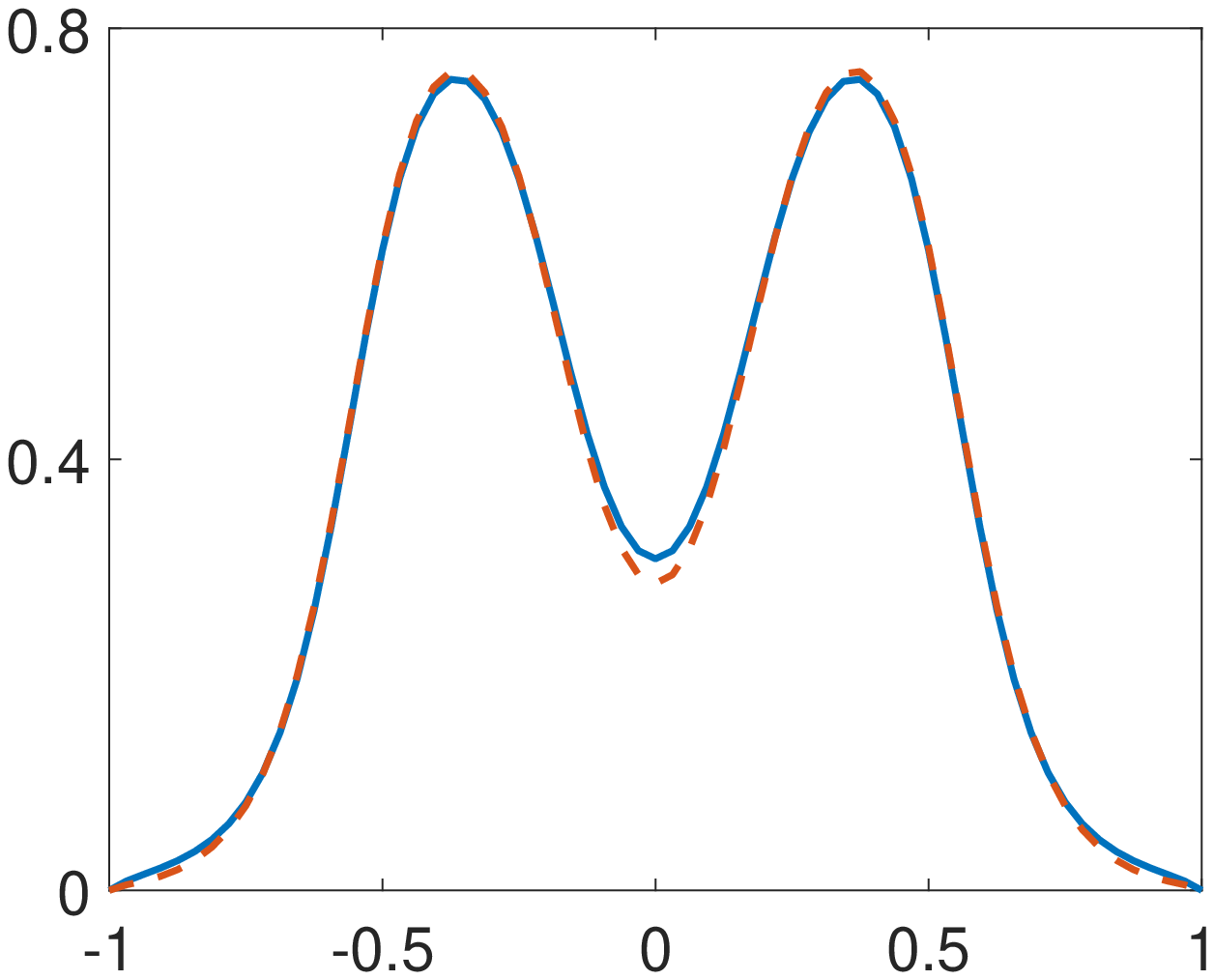}
}
	\caption{Plots of $u_0(x)$ in (a) and the approximated saddle points by SD (blue solid) and GPSD (red dashed) for the cases (i)--(iii) in (b)--(d), respectively.}
	\label{fig2}
\end{figure}

We employ SD and GPSD to compute saddle points and solution landscapes of the following nonlocal phase-field equation \cite{DuoWan,SonXu,YuZhe}
\begin{equation}\label{pf}
\partial_t u=F(u):=-\kappa\Big((-\Delta)^{\alpha/2}u+\frac{1}{\eta^2}(u-1)(u-1/2)u \Big) 
\end{equation}
on $x\in (-1,1)$, equipped with the initial and nonlocal boundary conditions
$$u(x,0)=u_0(x),~~x\in (-1,1);~~u(x,t)=0,~~x\in (-\infty,-1]\cup[1,\infty),~~t\geq 0.  $$
Here $\kappa$ represents the ``elastic'' relaxation time, $\eta$ is the interface parameter and the fractional Laplacian operator is defined as \cite{Eli,Lis,SheShe}
$$(-\Delta)^{\alpha/2}u(x,t):=C_\alpha\, \mbox{P.V.}\int_{-\infty}^\infty\frac{u(x,t)-u(y,t)}{\vert x-y\vert^{1+\alpha}}dy$$
where 
$$C_\alpha:=\frac{2^{\alpha-1}\alpha \Gamma(\alpha/2+1/2) }{\pi^{1/2}\Gamma(1-\alpha/2)}.$$
To discretize model (\ref{pf}), we adopt the numerical scheme of $(-\Delta)^{\alpha/2}$ in \cite{DuoZha} with the mesh size $h=2^{-5}$. Furthermore, we set $\alpha=1.5$, tol$_l$=0.05, tol$_u$=0.15, $N_m=20000$, tol$_x$=$1\times 10^{-5}$, $N_{sam}=N_{new}=120$, $u_0(x)=0.5(1-x^2)$ as shown in Figure \ref{fig2}(a), $\kappa=0.02$, $\Delta=0.01$ and the shrinking dimer with exponential decay $l_n=e^{-n\tau}l_0$. The initial values $\{v_{i,0}\}_{i=1}^k$ are chosen as the orthonormal eigenvectors of the Jacobian of $F$ at $u_0$. Appropriate smooth curves serve as the samples of the training set. We first take $\tau=l_0=0.00025$ to compute the saddle points of model (\ref{pf}) under different values of $1/\eta^2$ and different indexes as follows
$$(\mbox{i}) \,\frac{1}{\eta^2}=30, \,\mbox{index}=1;~(\mbox{ii}) \,\frac{1}{\eta^2}=80, \,\mbox{index}=2;~(\mbox{iii})\, \frac{1}{\eta^2}=120, \,\mbox{index}=3. $$  
Numerical results are presented in Figure \ref{fig2} and Table \ref{table4}, which lead to the same conclusions as in the first example. The we combine both SD and GPSD under $\tau=l_0=0.0005$ with the downward search algorithm proposed in \cite{YinPRL,YinSCM} to demonstrate the construction of the solution landscapes. We fix $1/\eta^2=120$ and $\kappa=0.02$ and start from an index-3 saddle point to construct the solution landscape of the nonlocal phase-field model (\ref{pf}). Numerical results are presented in Figure \ref{fig3}, which indicate that both methods generate similar solution landscapes, which demonstrates the effectiveness of the GPSD in solution landscape constructions.

	\begin{table}[h]
\setlength{\abovecaptionskip}{0pt}
\centering
\caption{Comparison of the number of force evaluation $N_f$ in SD and GPSD for computing the saddle points of the nonlocal phase field model (\ref{pf})}
\begin{tabular}{ccc} \cline{1-3}
&$N_f$ in SD&$N_f$ in GPSD  \\ \cline{1-3}		
(i)&18225&960	\\
(ii)& 95125&3480		\\
(iii)&117516 &4200	     \\
				\hline
			\end{tabular}
			\label{table4}
		\end{table}	

\begin{figure}[H]
	\setlength{\abovecaptionskip}{0pt}
	\centering	
	\includegraphics[width=4in,height=4in]{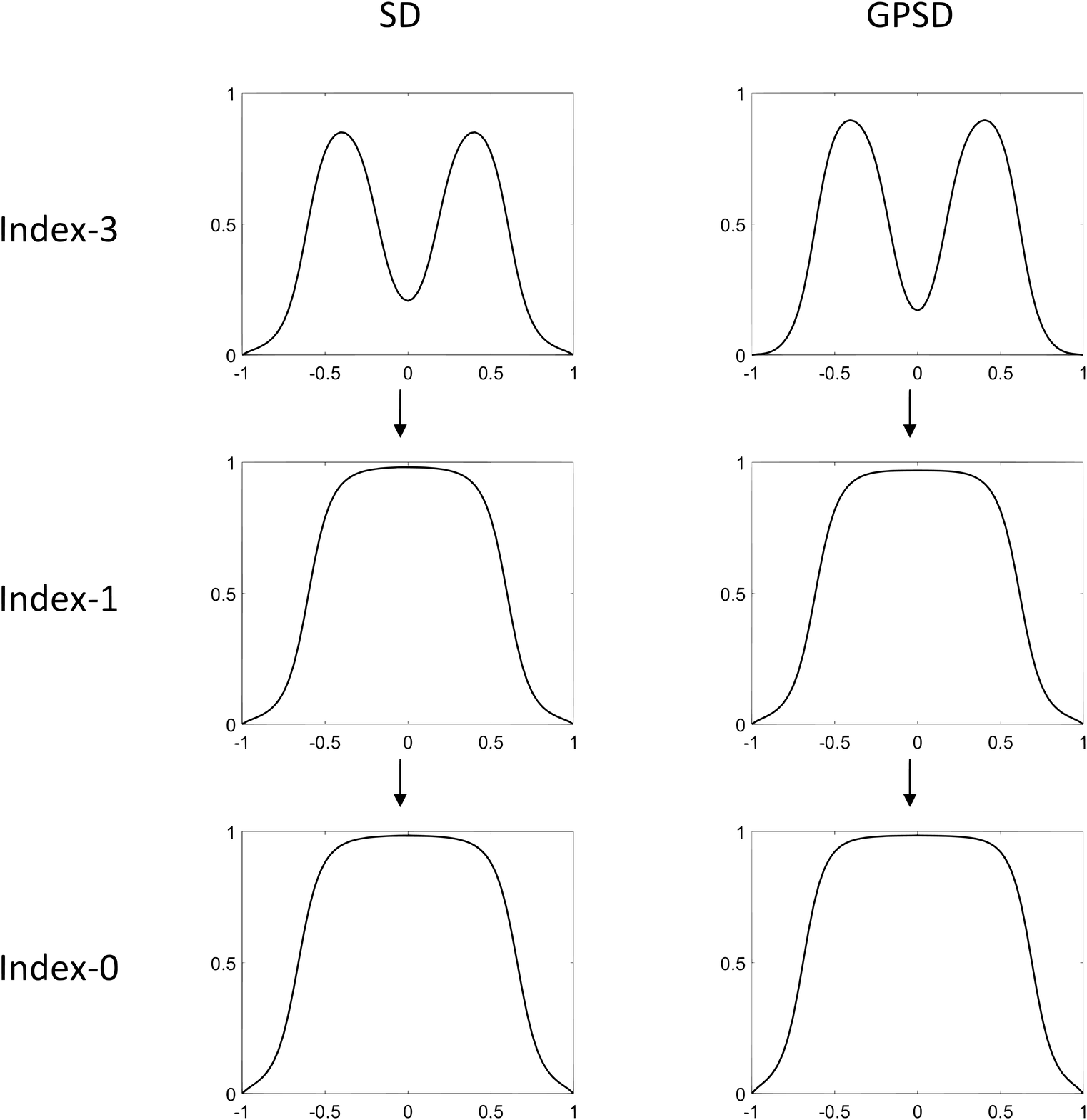}

	\caption{Solution landscape of nonlocal phase-field model (\ref{pf}) constructed by (left) SD and (right) GPSD based downward search algorithm.}
	\label{fig3}
\end{figure}

\section{Concluding remarks}
In this work we propose a model-free shrinking-dimer saddle dynamics and a corresponding sequential learning algorithm for finding any-index saddle points of complex models and constructing their solution landscapes. This data-driven approach avoids employing the exact form of the model and could significantly reduce the number of queries of the force that may be expensive or time-consuming, and the active learning mechanism along the latent trajectory increases the value of each sampling. More importantly, with the proposed algorithm, we can efficiently construct the model-free solution landscape, in which only the informative saddle points and minima are computed, so that one can avoid the sampling on the whole energy landscape of complex models. 

There are several potential extensions of the current method. For instance, different covariance kernels in the Gaussian process or  different learning methods such as the deep learning method could be applied to train the surrogate model. The proposed data-driven method could be further applied for constructing the solution landscapes of more complex problems, and how to utilize the structure of the solution landscape to further reduce the number of queries of the force is an interesting but challenging problem that will be investigated in the near future.

\bmhead{Acknowledgments}
We greatly appreciate Dr. Siwei Duo for providing the MATLAB code of computing the fractional Laplacian in \cite{DuoZha}.
\section*{Declarations}

\textbf{Funding}: This work was partially supported by the National Key R\&D Program of China No. 2021YFF1200500 and the National Natural Science Foundation of China No. 12225102, 12050002, 12288101. \\[0.1in]
\textbf{Conflict of interest}: The authors have no conflict of interest.\\[0.1in]
\textbf{Authors' contributions}: 
\textbf{Lei Zhang}: Conceptualization, Funding acquisition, Methodology, Supervision, Writing -- review \& editing. \textbf{Pingwen Zhang}: Funding acquisition, Supervision. 
\textbf{Xiangcheng Zheng}: Formal analysis, Funding acquisition, Investigation, Methodology, Visualization, Writing -- original draft.

\bigskip


\end{document}